\documentclass{article}
\usepackage[utf8]{inputenc}
\usepackage{hyperref}
\usepackage[round]{natbib}
\usepackage{amsmath,amssymb}

\begin{document}
\title{Martingale marginals do not always determine convergence}
\author{Jim Pitman \footnote{Department of Statistics, University of California, 367 Evans Hall, Berkeley, CA 94720-3860, U.S.A.} }
\date{\today}
\maketitle

\newcommand {\reals} {\mathbb{R}}
\newcommand {\complex} {\mathbb{C}}
\newcommand {\FF} {\mathcal{F}}
\newcommand {\GG} {\mathcal{G}}
\newcommand{\eq}{\begin{equation}}
\newcommand{\en}{\end{equation}}
\newcommand{\hf}{ \mbox{${\frac{1}{2}}$}}
\newcommand{\Ito}{It\^o}
\newcommand{\bV}{ {\bf V} }
\newcommand{\tilV}{ {   \tilde{V} }}
\newcommand{\br}{ \mbox{$\scriptstyle{\rm br}$}}
\newcommand{\re}[1]{\mbox{(\ref{#1})}}
\newcommand{\lb}[1]{\label{#1}}
\newcommand{\eval}{\left|\begin{array}{c}\\\\\end{array}\right.\!\!\!\!\!}
\def\qed{\mbox{\rule{0.5em}{0.5em}}}

\def\proof{\noindent{\bf Proof.\ \ }}
\def\note{\noindent{\bf Note.\ \ }}
\def\endpf{$\Box$}
\font\bb=msbm10
\def\bR{\hbox{\bb R}}
\newcommand {\fup}{ \mbox{$\phi^{\uparrow}(\alpha,$}}
\newcommand {\fdo}{ \mbox{$\phi^{\downarrow}(\alpha,$}}
\newcommand {\lamdel}{ \mbox{$\Lambda_{00}^{\delta}$}}

\newfont{\msbm}{msbm10 at 12pt}
\newfont{\eusb}{eusb10}
\newfont{\eusm}{eusm10}
\newfont{\eurb}{eurb10}
\newfont{\eurm}{eurm10}
\newfont{\eufb}{eufb10}
\newfont{\eufm}{eufm10}

\newcommand {\ints} {\mbox{\msbm\symbol{'132}}}
\newcommand {\IZ} {\mbox{\msbm\symbol{'132}}}
\newcommand {\rls} {\mbox{\msbm\symbol{'122}}}
\newcommand {\IR} {\mbox{\msbm\symbol{'122}}}
\newcommand {\IP} {\mbox{\msbm P}}
\newcommand {\IC} {\mbox{\msbm C}}
\newcommand {\rats} {\mbox{\msbm\symbol{'121}}}
\newcommand {\IQ} {\mbox{\msbm\symbol{'121}}}
\newcommand {\nats} {\mbox{\msbm\symbol{'116}}}
\newcommand {\IN} {\mbox{\msbm\symbol{'116}}}
\newcommand {\NN} {\mbox{\msbm\symbol{'116}}}
\newcommand {\IA}{ \mbox{$\cal A$}}
\newcommand{\con}{\rightarrow}
\newcommand{\te}{\rightarrow}
\newcommand{\ed}{\mbox{$ \ \stackrel{d}{=}$ }}
\newcommand{\giv}{\,|\,}
\newcommand{\convd}{\mbox{$ \ \stackrel{d}{\te}$ }}
\newcommand{\eps}{\varepsilon}
\newcommand {\nnn}{ {\bf n} }
\newcommand {\nnnn}{ \mbox{$(n_1, \ldots , n_k )$}}
\newcommand {\ppp}{ p({\bf n})}
\newcommand {\pppp} { \mbox{$p(n_1, \ldots , n_k )$}}
\newcommand {\sigman}{ \Sigma ({\bf n } ) }
\newcommand {\nseq}{ ( {\bf N }_1, \cdots, {\bf N } _ n ) }
\newcommand {\Nns}{ {\bf N } _ n  }

\newcommand{\EE}{{\cal E}}
\newcommand{\cd}{\mbox{$\ \stackrel{d}{\rightarrow} \ $}}
\newcommand{\cp}{\mbox{$\stackrel{p}{\rightarrow}$}}
\newcommand{\spten}{\,\,\,\,\,\,\,\,\,\,}

\section{Introduction}
Recent work of Marc Yor and coauthors 
\citeyearpar{yor:2011-5} 
has drawn attention to how properties of a martingale are related to its family of marginal distributions.
A fundamental result of this kind is Doob's martingale convergence theorem: 
\begin{itemize}
\item if the marginal distributions ($\mu_n, n \ge 0)$ of a discrete time martingale $(M_n, n \ge 0 )$ are such 
that $\int |x| \mu_n(dx)$ is bounded, then $M_n$ converges almost surely. 
\end{itemize}
Other well known results relating the behavior of a discrete time martingale $M_n$ to its marginal laws $\mu_n$ are:
\begin{itemize}
\item for each $p > 1$, the sequence $\int |x|^p \mu_n(dx)$ is bounded if and only if $M_n$ converges in $L^p$;
\item $\lim_{y \te \infty} \sup_{n} \int_{|x| > y } |x| \mu_n (dx) = 0$, that is $(M_n)_{n \ge 0}$ is uniformly integrable, 
if and only if $M_n$ converges in $L^1$.
\end{itemize}
We know also from L\'evy that if $\mu_n$ is the distribution of a partial sum $S_n$ of independent random variables,
and $\mu_n$ converges in distribution as $n \te \infty$, then $S_n$ converges almost surely.
These results can be found in most modern graduate textbooks in probability. See for instance Durrett \citeyearpar{durrett-probability2010}.

What  if the marginals of a martingale converge in distribution?  Does that imply the martingale converges a.s?
B\'aez-Duarte \citeyearpar{baez-duarte71} and Gilat \citeyearpar{gilat72} gave examples of martingales that converge in probability but not almost surely. 
So the answer is to this question is no.
But worse than that, there is a sequence of  martingale marginals converging in distribution, such that some martingales with these marginals
 converge almost surely, while others diverge almost surely. So by mixing, the probability of convergence of a martingale with these marginals can be any number in $[0,1]$.
Moreover, the same phenomemon can be exhibited for convergence in probability:
there is a sequence of  martingale marginals converging in distribution, such that some martingales with these marginals converge in probability, but others do not.

The purpose of this brief note is to record these examples, and to draw attention to the following problems which they raise:

\begin{itemize}
\item[1.] What is a necessary and sufficient condition on martingale marginals for every martingale with these marginals to converge almost surely?

\item[2.] What is a necessary and sufficient condition on martingale marginals for every martingale with these marginals to converge in probability?

\end{itemize}

Perhaps the condition for almost sure convergence is Doob's $L^1$-bounded condition.  But this does not seem at all obvious. 
Moreover, $L^1$-bounded is not the right condition for convergence in probability: convergence in distribution to a point mass is obviously sufficient, 
and this condition can hold for marginals that are not bounded in $L^1$.
See also Rao \citeyearpar{rao79} for treatment of some other problems related to non-$L^1$-bounded martingales.

\section{Examples}
\subsection{Almost sure convergence}

This construction extends and simplifies the  construction by Gilat 
\citeyearpar[\S 2]{gilat72} of a martingale which converges in probability but not almost surely, with increments in the set $\{-1,0,1\}$
See also B\'aez-Duarte \citeyearpar{baez-duarte71} for an earlier construction with unbounded increments, based on the double or nothing game 
instead of a random walk.

Let $(S_n, n = 0,1,2, \ldots )$ be a simple symmetric random walk started at $S_0 = 0$, with $(S_{n+1}-S_n, n = 0,1,2, \ldots)$ a sequence of
independent $U(\pm 1)$ random variables, where $U(\pm 1)$ is the uniform distribution on the set $\{\pm 1\}:= \{-1, +1\}$.
Let $0 = T_0 < T_1 < T_2 < \cdots$ be the successive times $n$ that $S_n = 0$. By recurrence of the simple random walk, $P(T_n < \infty ) = 1$ for every $n$.
For each $k = 1, 2 , \ldots $ let $M^{(k)}$ be the process which follows the walk $S_n$ on the random interval  $[T_{k-1}, T_k]$ of its $k$th excursion away from $0$, and is otherwise identically $0$:
$$
M_n^{(k)} := S_n 1( T_{k-1} \le n \le T_k )
$$
where $1(\cdots)$ is an indicator random variable with value $1$ if $\cdots$ and $0$ otherwise. Each of these processes $M^{(k)}$ is a martingale relative to the filtration $(\FF_n)$ generated by
the walk $(S_n)$, by Doob's optional sampling theorem.
Now let $(A_k)$ be a sequence of events such that the $\sigma$-field $\GG_0$ generated by these events is independent of the walk $(S_n, n \ge 0)$, 
and set
$$
M_n:= \sum_{k=1}^\infty M_n^{(k)} 1(A_k) 
$$
So $M_n$ follows the path of $S_n$ on its $k$th excursion away from $0$ if $A_k$ occurs, and otherwise $M_n$ is identically $0$.
Let $\GG_n$ for $n \ge 0$ be the $\sigma$-field generated by $\GG_0$ and $\FF_n$. Then it is clear that $(M_n, \GG_n)$ is a martingale, no matter what 
choice of the sequence of events $(A_k)$ independent of $(S_n)$. The distribution of $M_n$ is determined by the formula
$$
P(M_n = x) = \sum_{k=1}^\infty P( S_n = x , T_{k-1} \le n \le T_k ) P(A_k)
$$
for all integers $x \ne 0$.
A family of martingales with the same marginals is thus obtained by varying the structure of dependence between the events $A_k$ for a given sequence of probabilities $P(A_k)$.
The only way that a path of $M_n$ can converge is if $M_n$ is eventually absorbed in state $0$. 
So if $N:=  \sum_{k} 1(A_k)$ denotes the number of events $A_k$ that occur,
$$
P( M_n \mbox { converges} ) = P( N < \infty ).
$$
Now take $P(A_k) = p_k$ for a decreasing sequence $p_k$ with limit $0$ but $\sum_k p_k = \infty$,  for instance $p_k = 1/k$.
Then $(A_k)$ can be constructed so that the $A_k$ are mutually indendent, and $P(N = \infty) = 1$ by the Borel-Cantelli lemma.
Or these events can be nested: 
$$A_1 \supseteq A_2 \supseteq A_3 \cdots $$ 
in which case 
$$
P(N \ge k ) = P(A_k) \downarrow 0 \mbox{ as $k \te \infty$},
$$
so $P(N= \infty)  = 0$ in this case.
Thus we obtain a sequence of marginal distributions for a martingale, such that some martingales with these marginals converge almost surely, while others
diverge almost surely.

\subsection{Convergence in probability}

Let us construct a martingale $M_n$ which converges in distribution, but not in probability,
following indications of such a construction by Gilat \citeyearpar[\S 1]{gilat72}.

This will be an inhomogeneous Markov chain with integer values, starting from $M_0 = 0$. 
Its first step will be to $M_1$ with $U(\pm 1)$ distribution. Thereafter, the idea is to force $M_n$ to alternate between the values $\pm 1$, with probability increasing to $1$
as $n \te \infty$. This achieves $U(\pm 1)$ as its limit in distribution, while preventing convergence in probability by the alternation.
The transition probabilities of $M_n$ are as follows:

\begin{align}
P( M_{n+1} = M_n \pm 1 \giv M_n \mbox { with } M_n \notin \{\pm 1 \} ) & = 1/2 \\
P( M_{n+1} = - 1  \giv M_n = 1 ) & = 1 - 2^{-n} \\
P( M_{n+1} = 2^{n+1} - 1  \giv M_n = 1 ) & = 2^{-n} \\
P( M_{n+1} = + 1  \giv M_n = -1 ) & = 1 - 2^{-n} \\
P( M_{n+1} = - 2^{n+1} + 1  \giv M_n = -1 ) & = 2^{-n}
\end{align}
The first line indicates that whenever $M_{n}$ is away from the two point set $\{\pm 1 \}$, it moves according to a simple symmetric random
walk, until it eventually gets back to $\{\pm 1 \}$ with probability one. Once it is back in $\{\pm 1 \}$, it is forced to alternate between these values,
with probability $1-2^{-n}$ for an alternation at step $n$, compensated by moving to $\pm ( 2^{n+1} - 1 )$ with probability $2^{-n}$. Since the probabilities $2^{-n}$ are summable, the Borel-Cantelli Lemma
ensures that with probability one only finitely many exits from $\{\pm 1 \}$ ever occur. After the last of these exits, the martingale eventually returns to $\{\pm 1\}$
with probability one. From that time onwards, the martingale flips back and forth deterministically between $\{\pm 1 \}$.

A slight modification of these transition probabilities, gives another martingale with the same marginal distributions which converges almost surely and hence in probability. With $M_0 = 0$ as before, the modified scheme is as follows:
\begin{align}
P( M_{n+1} = M_n \pm 1 \giv M_n \mbox { with } M_n \notin \{\pm 1 \} ) & = 1/2 \\
P( M_{n+1} = 1  \giv M_n = 1 ) & = 1 - 2^{-n} \\
P( M_{n+1} = 2^{n+1} - 1  \giv M_n = 1 ) & = 2^{-n} p_n \\
P( M_{n+1} = - 2^{n+1} + 1  \giv M_n = 1 ) & = 2^{-n} q_n \\
P( M_{n+1} = - 1  \giv M_n = -1 ) & = 1 - 2^{-n} \\
P( M_{n+1} = - 2^{n+1} + 1  \giv M_n = -1 ) & = 2^{-n} p_n \\
P( M_{n+1} = 2^{n+1} - 1  \giv M_n = -1 ) & = 2^{-n} q_n 
\end{align}
where
$$
p_n := 1/ ( 2 - 2^{-n} ) \mbox { and } q_n := 1 - p_n
$$
are chosen so that the distribution with probability $p_n$ at $2^{n+1} - 1$ and $q_n$ at $-2^{n+1} + 1$ has mean 
$$
p_n ( 2^{n+1} - 1 )  + q_n ( -2^{n+1} + 1 ) = 1.
$$
In this modified process, the alternating transition out of states $\pm 1$ is replaced by holding in these states, while the previous compensating
moves to $\pm ( 2^{n+1} - 1 ) $ are replaced by a nearly symmetric transitions from $\pm 1$ to these values. This preserves the martingale property,
and also preserves the marginal laws. But the previous argument for eventual alternation now shows that the modified martingale is eventually absorbed almost surely
in one of the states $\pm 1$. So the modified martingale converges almost surely to a limit which has $U(\pm 1 )$ distribution.

These martingales $(M_n)$ have jumps that are unbounded. Gilat \citeyearpar[\S 2]{gilat72} left open the question of whether there exist martingales with uniformly bounded increments which converge in distribution but not in probability. But such martingales can be created by a variation of the first construction of $(M_n)$ above, as follows.

Run a simple symmetric random walk starting from $0$. Each time the random walk makes an alternation between the two states $\pm 1$, make the 
walk delay for a random number of steps in its current state in $\pm 1$ before continuing, for some rapidly increasing sequence of random delays. 
Call the resulting martingale $M_n$. 
So by construction, $M_1$ has $U(\pm)$ distribution, 
$$M_n = (-1)^{k-1} M_1 \mbox{  for } S_k \le n \le T_k$$ 
for some increasing sequence of randomized stopping times 
$$1 = S_1 < T_1 < S_2 < T_2 < \cdots,$$
and during the $k$th {crossing interval}  $[T_k, S_{k+1}]$
the process $M_n$ follows a simple random walk path starting  in state $(-1)^{k-1} M_1$ and stopping when it first reaches state $(-1)^{k} M_1$.

The claim is that a suitable construction of the delays $T_k - S_k$ will ensure that the distribution of $M_n$ converges to $U(\pm1)$, while there is almost deterministic
alternation for large $k$ of the state $M_{t_k}$ for some rapidly increasing deterministic sequence $t_k$.
To achieve this end, let $t_1 = 1$ and suppose inductively for $k = 1, 2, \ldots$ that $t_k$ has been chosen
so that 
\begin{equation}
\label{alt}
P( M_{t_k} = (-1)^{k-1}M_1) > 1 - \epsilon_k \mbox{ for some $\epsilon_k \downarrow 0$ as $k \te \infty$}.
\end{equation}
Here $M_1 \in \{\pm 1 \}$ is the first step of the simple random walk.
The random number of steps required for random walk crossing between states $\pm 1$ is a.s. finite.
So having defined $t_k$, we can choose  an even integer $t_{k+1}$ so large, that $t_{k+1}/2 > t_k$ and all of the following events occur with probability at least  $1 - \epsilon_{k+1}$:
\begin{itemize}
\item $M_{t_{k+1}/2 } = (-1)^{k-1} M_1$, meaning that the $(k-1)$th crossing between $\pm 1$ has been completed by time $S_{k} < t_{k+1}/2$;
\item the $k$th crossing is started at time $T_k$ that is uniform on $[t_{k+1}/2, t_{k+1})$ given $S_{k} < t_{k+1}/2$;
\item the $k$th crossing is completed at time $S_{k+1} < t_{k+1}$, so $M_n = (-1)^k M_1$ for $S_{k+1} \le n \le t_{k+1}$.
\end{itemize}
Moreover, $t_{k+1}$ can be chosen so large that the uniform random start time of the $k$th crossing given $S_{k} < t_{k+1}/2$
ensures that  also
$$
P( M_n \in \{\pm1\} ) \ge  1 - 2 \epsilon_k \mbox{   for all } t_k \le n \le t_{k+1}
$$
because with high probability the length $S_{k+1} - T_k$ of the $k$th crossing 
is negligible in comparison with the length $t_{k+1}/2$ of the interval $[t_{k+1}/2 , t_{k+1}]$  in which this crossing is arranged to occur.
It follows from this construction that $M_n$ converges in distribution to $U(\pm1 )$,
while the forced alternation \eqref{alt} prevents $M_n$ from having a limit in probability.

A feature of the previous example is that $\sup_n M_n = - \inf_n M_n = \infty$ almost surely, since in the end every step 
of the underlying simple symmetric random walk is made by the time-changed martingale $M_n$.
A similar example can be created from a standard Brownian motion $(B_t, t \ge 0)$ 
using a predictable $\{0,1\}$-valued process $(H_t, t \ge 0)$ 
to create successive switching between and holding in states $\pm 1$ so that the martingale
$$
M_t := \int_0^t H_t dB_t
$$
converges in distribution to $U(\pm 1)$ while not converging in probability. In this example, $\int_0^\infty H_t dt = \sup_t M_t = - \inf_t M_t = \infty$ almost surely.

\paragraph{Acknowledgement}
Thanks to David Aldous for drawing my attention to Gilat \citeyearpar{gilat72}.

\def\cprime{$'$} \def\polhk#1{\setbox0=\hbox{#1}{\ooalign{\hidewidth
  \lower1.5ex\hbox{`}\hidewidth\crcr\unhbox0}}} \def\cprime{$'$}
  \def\cprime{$'$} \def\cprime{$'$}
  \def\polhk#1{\setbox0=\hbox{#1}{\ooalign{\hidewidth
  \lower1.5ex\hbox{`}\hidewidth\crcr\unhbox0}}} \def\cprime{$'$}
  \def\cprime{$'$} \def\polhk#1{\setbox0=\hbox{#1}{\ooalign{\hidewidth
  \lower1.5ex\hbox{`}\hidewidth\crcr\unhbox0}}} \def\cprime{$'$}
  \def\cprime{$'$} \def\cydot{\leavevmode\raise.4ex\hbox{.}} \def\cprime{$'$}
  \def\cprime{$'$} \def\cprime{$'$} \def\cprime{$'$} \def\cprime{$'$}
  \def\cprime{$'$} \def\cprime{$'$} \def\cprime{$'$} \def\cprime{$'$}
  \def\cprime{$'$}

\end{document}